# Incorporating a contrast in the Bayesian formula: What consequences for the MAP estimator and the posterior distribution? Applications in spatial statistics

S. Soubeyrand*†‡, F. Carpentier*, N. Desassis§ and J. Chadœuf*

November 15, 2018

**Abstract.** In order to estimate model parameters and circumvent possible difficulties encountered with the likelihood function, we propose to replace the likelihood in the formula of the posterior distribution by a function depending on a contrast. The properties of the contrast-based (CB) posterior distribution and MAP estimator are studied to understand what the consequences of incorporating a contrast in the Bayesian formula are. We show that the proposed method can be used to make frequentist inference and allows the reduction of analytical calculations to get the limit variance matrix of the estimator. For specific contrasts, the CB–posterior distribution directly approximates the limit distribution of the estimator; the calculation of the limit variance matrix is then avoided. Moreover, for these contrasts, the CB–posterior distribution can also be used to

*INRA, UR546 Biostatistique et Processus Spatiaux, Domaine St Paul, F-84914 Avignon, France

†INRA, Agro ParisTech, UMR1290 BIOGER-CPP, BP01 F-78850 Thiverval-Grignon, France

‡Corresponding author: samuel.soubeyrand@avignon.inra.fr

§MERE project, INRIA/INRA, UMR INRA/SupAgro Analyse des Systèmes et Biométrie, 2 place Viala, F-34060 Montpellier, France

make inference in the Bayesian way. The method is applied to three spatial data sets.

**Key words.** Frequentist estimation; Quasi-Bayesian estimation; Spatial models.

# 1 Introduction

In both the frequentist and the Bayesian viewpoints, the likelihood function has become the major component of statistical inference under a parametric model. Its use, however, has drawbacks in specific situations. First, it may be impossible to write down the likelihood in a numerically tractable form; see the cases of Boolean models (Van Lieshout and Van Zweit, 2001), Markov point processes (Møller, 2003), Markov spatial processes (Guyon, 1985) and spatial generalized linear mixed models (spatial GLMM; Diggle et al., 1998) where multiple integrals cannot be reduced due to spatial dependences. Second, the likelihood may not be completely appropriate because of the associated assumptions. For instance, the likelihood is built under an assumption on the distribution of data, but such an assumption may be tricky to specify in case of insufficient information as in classical geostatistics (Chilès and Delfiner, 1999); see also McCullagh and Nelder (1989, chap. 9). In the same vein, every data are assumed to have the same weights in the likelihood, but the influence of outliers may be too large according to the analyst (Markatou, 2000).

The difficulties encountered with the likelihood can be circumvented with existing Bayesian and frequentist procedures.

- There are procedures which use conditional simulation to numerically approximate the likelihood. For instance, the Markov chain Monte Carlo algorithm (MCMC; Robert and Casella, 1999), for example, allows the approximation of the posterior distribution for Markov point processes (Møller, 2003) and spatial GLMMs (Diggle et al., 1998). The Markov chain expectation maximization algorithm (MCEM; Wei and Tanner, 1990) allows the maximization of the likelihood for Boolean models (Van Lieshout and Van Zweit, 2001) and spatial GLMMs (Zhang, 2002).

- There are procedures where the likelihood function is simplified or replaced. For example, the pseudo-likelihood, which only accounts for local dependence structures, is used instead of the likelihood for Markov point processes (Møller, 2003) and Markov spatial processes (Besag, 1975; Guyon, 1985). The generalized least squares estimation, which does not rely on assumptions on the distribution of data, is used in geostatistics; see Chilès and Delfiner (1999, chap. 2-3) and Stein (1999, chap. 1). Other procedures belonging to this category are: the weighted likelihood maximization (Markatou, 2000), the method of moments, the M-estimation (Serfling, 2002), the approximate Bayesian computation (ABC; Beaumont et al., 2002), the quasi-likelihood maximization (McCullagh and Nelder, 1989) and the quasi-Bayesian likelihood method (Lin, 2006).

In the quasi-Bayesian likelihood approach, the likelihood appearing in the posterior distribution formula is replaced by a quasi-likelihood which does not rely on distribution assumptions. Then, the posterior distribution which is obtained is used to make inference as in classical Bayesian situations. In this communication we propose to generalize this approach: the likelihood in the posterior distribution formula is replaced by a function of a contrast.

A contrast is a function of the model parameters and the observed data which is minimized to estimate the parameters (Dacunha-Castelle and Duflo, 1982). The minimum contrast approach is a generic estimation method which was developed in a frequentist perspective. The maximum likelihood estimation as well as the maximum pseudo, weighted or quasi likelihood estimation, the diverse least squares methods, the method of moments and the M-estimation can be formulated as minimum contrast estimation problems.

Thus, the procedure which is proposed —replacing the likelihood by a function of a contrast in the Bayesian formula— includes the classical Bayesian approach (here and thereafter, "classical" refers to "likelihood-based") and the quasi-Bayesian approach of Lin (2006). This procedure provides a contrast-based (CB) posterior distribution which does not coincide, in the general case, with the classical posterior distribution. In this paper, we investigate what are the posterior distribution and the MAP (maximum a posteriori) estimator based on a contrast.

Under mild conditions on the prior distribution, we show that the CB–MAP estimator inherits the asymptotic properties (consistency and asymptotic normality) of the minimum contrast estimator, as the classical MAP estimator inherits the asymptotic properties of the maximum likelihood estimator (Caillot and Martin, 1972). The limit variance matrix of the normalized estimator is $I_\theta^{-1}\Gamma_\theta I_\theta^{-1}$ where $\Gamma_\theta$ is the limit variance of the gradient of the contrast and $I_\theta$ is the limit Hessian matrix of the contrast.

Moreover, we show that the CB–posterior distribution is asymptotically equivalent to a normal distribution whose variance matrix is $I_\theta^{-1}$. Therefore, when building the contrast, particular attention must be paid to satisfy, if possible, $I_\theta^{-1}\Gamma_\theta I_\theta^{-1} = I_\theta^{-1}$. Indeed, with such a contrast, inference can be made without computing matrices $\Gamma_\theta$ and $I_\theta$: the posterior distribution can either be used as a limit distribution in a frequentist viewpoint or be used to make inference in the Bayesian way. When building a contrast satisfying $I_\theta^{-1}\Gamma_\theta I_\theta^{-1} = I_\theta^{-1}$ is not possible, the CB–posterior distribution can nevertheless be used to estimate $I_\theta^{-1}$. Thus, the computation of the limit Hessian matrix of the contrast is avoided.

To summarize, the present study shows the consequences of replacing the likelihood by a function of a contrast. It also provides an estimation method which has advantages over existing methods exploited to circumvent difficulties encountered with the likelihood. First, it does not require a simulation–based algorithm as the MCMC, MCEM or ABC algorithms. Second, it inherits the richness of the minimum contrast approach (there are many types of contrast: likelihood, least squares, moments...). Third, compared to the classical contrast method, the computation of the derivatives of the contrast is limited. Fourth, when $I_\theta^{-1}\Gamma_\theta I_\theta^{-1} = I_\theta^{-1}$, the CB-posterior distribution can be directly used to make inference either in the frequentist perspective or in the Bayesian perspective. However, the method which is proposed has also drawbacks. In particular, building a contrast which exploits a large part of the information in the data, as the likelihood does, is not obvious. Besides, building a contrast satisfying $I_\theta^{-1}\Gamma_\theta I_\theta^{-1} = I_\theta^{-1}$ asks analytical work which can be time consuming. Furthermore, obtaining such a contrast is not always possible.

The article is organized as follows. The classical minimum contrast method of estimation is recalled in section 2 and examples are given. The method that we propose is presented in section 3, and its properties are derived. Then, the method is applied in section 4 to simulated and real cases dealing with spatial statistics (estimation of the range parameter of a variogram; estimation of the parameters of a Markovian spatial process; and estimation of the parameters of an autosimilar model used to describe soil roughness). The three cases illustrate the application of the method when the parameter has one or several components and when $I_\theta^{-1}\Gamma_\theta I_\theta^{-1}$ is equal to or different from $I_\theta^{-1}$.

# 2 Recall: Classical minimum contrast estimation

## 2.1 Estimator and asymptotic properties

Detailed information on minimum contrast estimation can be found in Dacunha-Castelle and Duflo (1982). Here, we avoid the complete notations. Consider a family of parametric models $\{P_\alpha : \alpha \in \Theta\}$ and samples of increasing sizes $t \in T \subset \mathbb{N}$, drawn from $P_\theta$. A contrast for $\theta$ is a random function $\alpha \mapsto U_t(\alpha)$ defined over $\Theta$, depending on a sample of size $t$, and such that $\{U_t(\alpha)\}_t$ converges in probability, as $t \to \infty$, to a function $\alpha \mapsto K(\alpha, \theta)$ which has a strict minimum at $\alpha = \theta$. The minimum contrast estimator is

$$\hat{\theta}_t = \mathrm{argmin}\{U_t(\alpha), \alpha \in \Theta\}.$$

Let us make the following classical assumptions:

$H_1$ : $\Theta \subset \mathbb{R}^p$, $p < \infty$, is compact and $\theta$ is in the interior of $\Theta$,

$H_2$ : $\alpha \mapsto K(\alpha, \theta)$ has a strict minimum at $\theta$,

$H_3$ : $\alpha \mapsto U_t(\alpha)$ is $C^2$ (it has two continuous derivatives) over $\Theta$,

$H_4$ : the normalized gradient vector $\sqrt{t}\mathbf{grad}U_t(\theta)$ (first derivatives of $U_t(\theta)$ with respect to $\theta$) converges in law to the normal distribution $\mathcal{N}(0, \Gamma_\theta)$:

$$\sqrt{t}\mathbf{grad}U_t(\theta) \to \mathcal{N}(0, \Gamma_\theta) \quad \text{in law as } t \to \infty,$$

$H_5$ : the Hessian matrix $\mathbf{H}U_t(\theta)$ (second derivatives of $U_t(\theta)$ with respect to $\theta$) converges in probability to an invertible matrix $I_\theta$:

$$\mathbf{H}U_t(\theta) \to I_\theta \quad \text{in probability as } t \to \infty,$$

$H_6$ : $\sup\limits_{||\beta||<\epsilon} |\mathbf{H}_{kl}U_t(\theta + \beta) - \mathbf{H}_{kl}U_t(\theta)| \to 0$ in probability, where $\epsilon > 0$ and $\mathbf{H}_{kl}$ is the component $(k, l)$, $1 \leq k, l \leq p$, of the Hessian operator.

Under these assumptions, the minimum contrast estimator is consistent and asymptotically normal: as $t \to \infty$,

- $\hat{\theta}_t$ converges in probability to $\theta$ and

- $\sqrt{t}(\hat{\theta}_t - \theta)$ converges in law to the Gaussian distribution $\mathcal{N}\left(0, I_\theta^{-1} \Gamma_\theta I_\theta^{-1}\right)$.

## 2.2 Examples

**Maximum likelihood.** Consider an i.i.d. sample $(X_i)_{1 \leq i \leq n}$ (here, $T = \mathbb{N}$), each element being drawn from the density $p_\theta(.)$. The likelihood function is $L_n = \prod_{0 \leq i \leq n} p_\theta(X_i)$ and the corresponding contrast is

$$U_n(\alpha) = -\frac{1}{n} \sum_{i \leq n} \log p_\alpha(X_i).$$

The limit function $K$ is the opposite of the Kullback information: $K(\alpha, \theta) = -E_\theta\{\log p_\alpha(X_i)\}$, the matrices $I_\theta$ et $\Gamma_\theta$ satisfy

$$I_\theta = \Gamma_\theta = E_\theta[\,\mathrm{grad}_\theta\{\log p_\theta(X_i)\}\,\mathrm{grad}_\theta\{\log p_\theta(X_i)\}'\,],$$

and the convergence in law simplifies into $\sqrt{n}(\hat{\theta}_n - \theta) \to \mathcal{N}\left(0, I_\theta^{-1}\right)$.

**Least squares.** Here we present the least-square method as a contrast method in the case of the estimation of a variogram. This case will be used as an illustration in the application section.

Consider a stationary Gaussian random field $X$ over $\mathbb{Z}^2$ with mean value zero and with parametric variogram $\gamma_\theta(h) = E_\theta\{(X_i - X_j)^2\}$, where $h = d(i, j)$ is the distance between $X_i$ and $X_j$ (Chilès and Delfiner, 1999). Assume that the sample is made on a square grid $\{i = (i_1, i_2) : 0 \leq i_1, i_2 \leq n\}$ with size $n^2$; the sample is denoted by $(X_i)_{0 \leq i_1, i_2 \leq n}$ where $i = (i_1, i_2)$ (here $T = \{n^2 : n \in \mathbb{N}\}$). The variogram can be estimated with the least square method (Chilès and Delfiner, 1999). In practice, the sample variogram $\hat{\gamma}$ is computed at each possible distance $h_l$ ($l \leq k$) between points: $\hat{\gamma}(h_l) = \frac{1}{2n_l} \sum_{(i,j) \in \mathcal{C}_l} (X_i - X_j)^2$, where $\mathcal{C}_l$ is the set of pairs of points separated by $h_l$ and $n_l = \#\mathcal{C}_l$, and the contrast between the sample variogram and the theoretical variogram

$$U_{n^2}(\alpha) = \frac{1}{2} \sum_{l \leq k} \{\hat{\gamma}(h_l) - \gamma_\alpha(h_l)\}^2 \tag{1}$$

is minimized. The limit function $K$ of the contrast is $K(\alpha, \theta) = \frac{1}{2} \sum_{l \leq k} \{\gamma_\theta(h_l) - \gamma_\alpha(h_l)\}^2$. In this context, the sample variogram $\{\hat{\gamma}(h_l)\}_{l \leq k}$ is unbiased with mean $\mu_\theta = \{\gamma_\theta(h_l)\}_{l \leq k}$ and $n\{\hat{\gamma}(h_l) - \gamma_\theta(h_l)\}_{l \leq k}$ is asymptotically normal with variance matrix denoted by $\Sigma_\theta$. It follows that $n(\hat{\theta}_n - \theta) \to \mathcal{N}(0, I_\theta^{-1} \Gamma_\theta I_\theta^{-1})$ where the component $(i,j)$ of $\Gamma_\theta$ is $\frac{\partial \mu_\theta'}{\partial \theta_i} \Sigma_\theta \frac{\partial \mu_\theta}{\partial \theta_j}$, the component $(i,j)$ of $I_\theta$ is $-\frac{\partial \mu_\theta'}{\partial \theta_i} \frac{\partial \mu_\theta}{\partial \theta_j}$ and $\mu_\theta'$ is the transpose of $\mu_\theta$.

**Pseudo-likelihood.** Here we present the pseudo-likelihood method as a contrast method in the case of the estimation of the parameters of a Markov random field. This case will be used as an illustration in the application section.

Consider a stationary Markov random field $X$ over $\mathbb{Z}^2$ with state space $\{0, 1\}$. Assume that the conditional probability of $X_i$ given $X_j, j \neq i$, satisfies

$$\begin{aligned} P_\theta(X_i \mid X_j, j \neq i) &= P_\theta(X_i \mid X_j, j \in V(i)) \\ &= \frac{\exp\left(\theta_1 X_i + \theta_2 \sum_{j \in V(i)} X_i X_j\right)}{\left\{1 + \exp\left(\theta_1 + \theta_2 \sum_{j \in V(i)} X_j\right)\right\}}, \end{aligned}$$

where $\theta = (\theta_1, \theta_2)$ is a pair of parameters and $V(i)$ is the set of the four nearest neighbors of $i$ (Guyon, 1985). We assume in the following that the Markov field is $\alpha$-mixing; this is satisfied if $\mid \theta_2 \mid \leq 1$ for example. Moreover, the field is observed on the square grid $\mathcal{I} = \{i = (i_1, i_2) : 0 \leq i_1, i_2 \leq n\}$ with size $n^2$ (here $T = \{n^2 : n \in \mathbb{N}\}$). The likelihood cannot be analytically calculated. Therefore, a pseudo-likelihood was proposed to make the inference (Guyon, 1985). The pseudo-likelihood is the product of the conditional probabilities $\prod_{i \in \mathcal{I}} P_\theta(X_i \mid X_j, j \neq i)$. The corresponding contrast is

$$U_{n^2}(\alpha) = -\frac{1}{n^2} \sum_{i \in \mathcal{I}} \log P_\alpha(X_i \mid X_j, j \in V(i)). \tag{2}$$

Let $\mathcal{W}$ denote the set of possible states for the neighborhood of any point 0, then the limit function of the contrast is

$$K(\alpha, \theta) = -\sum_{w \in \mathcal{W}} \sum_{x \in \{0,1\}} \log P_\alpha\{x \mid X_i = w_i, i \in V(0)\} P_\theta\{x \mid X_i = w_i, i \in V(0)\} P_\theta(w).$$

Moreover, $n(\hat{\theta}_n - \theta) \to \mathcal{N}(0, I_\theta^{-1} \Gamma_\theta I_\theta^{-1})$ where $I_\theta = \text{var}(Z_0)$, $\Gamma_\theta = M_0 + 4 \sum_{0 \leq i_1, i_2 \leq 2} M_i$, $M_i = \text{cov}(Z_0, Z_i)$, $i \in \mathcal{I}$, and vectors $Z_i$ satisfy

$$Z_i = \left( X_i - \frac{\exp\left(\theta_1 + \theta_2 \sum_{j \in V(i)} X_j\right)}{1 + \exp\left(\theta_1 + \theta_2 \sum_{j \in V(i)} X_j\right)} \right) \begin{pmatrix} 1 \\ \sum_{j \in V(i)} X_j \end{pmatrix}.$$

# 3 Incorporating a contrast in the Bayesian formula

## 3.1 Posterior distribution and MAP estimator based on a contrast

In the Bayesian framework, a prior distribution denoted $c(\cdot)$ is defined over $\Theta$. Let $(X_i)_{i \leq t}$ be a sample of size $t$ drawn from the distribution $P_\theta$, then the posterior distribution is

$$p(\theta \mid X_i, i \leq t) = \frac{P_\theta(X_i, i \leq t) c(\theta)}{\int_\Theta P_\alpha(X_i, i \leq t) c(\alpha) d\alpha}$$
$$= \frac{\exp(-t U_t(\theta)) c(\theta)}{\int_\Theta \exp(-t U_t(\alpha)) c(\alpha) d\alpha}$$

where $P_\theta(X_i, i \leq t)$ is the likelihood and $U_t(\alpha) = -\frac{1}{t} \log P_\alpha(X_i, i \leq t)$ is the corresponding contrast (see the first example presented above).

For the estimation of $\theta$, we propose to replace the contrast associated with the likelihood in the Bayesian formula written above by any contrast. We obtain a contrast-based (CB) posterior distribution denoted $p_t(\alpha)$:

$$p_t(\alpha) = \frac{\exp(-t U_t(\alpha)) c(\alpha)}{\int_\Theta \exp(-t U_t(\beta)) c(\beta) d\beta}. \tag{3}$$

The CB–MAP estimator obtained by maximizing $p_t(\cdot)$ is denoted

$$\tilde{\theta}_t = \mathrm{argmax}\{p_t(\alpha), \alpha \in \Theta\}.$$

$\tilde{\theta}_t$ is at the minimum of $\alpha \mapsto U_t(\alpha) - (1/t)\log c(\alpha)$, and does not coincide in the general case with the classical minimum contrast estimator $\hat{\theta}_t = \mathrm{argmin}\{U_t(\alpha), \alpha \in \Theta\}$.

In what follows we investigate the behavior of the CB–MAP estimator and the CB–posterior distribution.

## 3.2 Consistency and asymptotic normality of the CB–MAP estimator

We noted above that the CB–MAP estimator $\tilde{\theta}_t$ is at the minimum of $\alpha \mapsto U_t(\alpha) - (1/t)\log c(\alpha)$. This function satisfies the definition of a contrast. Consequently, convergence properties of $\tilde{\theta}_t$ can be easily obtained by using the contrast theory. Assume that the hypotheses listed in section 2 are satisfied. Let us assume in addition that the prior distribution $c(\cdot)$ is differentiable and strictly positive over $\Theta$. It can be shown that, as $t \to \infty$,

- $\tilde{\theta}_t$ converges in probability to $\theta$ and

- $\sqrt{t}(\tilde{\theta}_t - \theta)$ converges in law to the Gaussian distribution $\mathcal{N}\left(0, I_\theta^{-1}\Gamma_\theta I_\theta^{-1}\right)$, where $I_\theta$ and $\Gamma_\theta$ are the matrices which were introduced when the classical minimum contrast method was presented:

$$\mathbf{H}U_t(\theta) \to I_\theta \quad \text{in probability as } t \to \infty$$
$$\sqrt{t}\mathbf{grad}U_t(\theta) \to \mathcal{N}(0, \Gamma_\theta) \quad \text{in law.}$$

## 3.3 Asymptotic deviation between $\tilde{\theta}_t$ and $\hat{\theta}_t$

The asymptotic deviation between the classical minimum contrast estimator $\hat{\theta}_t$ and the CB–MAP estimator $\tilde{\theta}_t$ is given by

$$\tilde{\theta}_t - \hat{\theta}_t = \frac{1 + o_{\text{proba}}(1)}{tc(\theta)} I_\theta^{-1} \mathbf{grad}c(\theta) \qquad (4)$$
$$= O_{\text{proba}}(t^{-1})\mathbf{1}_p.$$

where $\mathbf{1}_p$ is the unit vector of size $p$ (the dimension of $\Theta$). Thus, the deviation between the two estimators is of order $1/t$.

Proof of (4). As $\tilde{\theta}_t$ satisfies $\mathbf{grad}p_t(\tilde{\theta}_t) = 0$,

$$0 = -tc(\tilde{\theta}_t)\mathbf{grad}U_t(\tilde{\theta}_t) + \mathbf{grad}c(\tilde{\theta}_t).$$

Then, applying a first order Taylor's expansion for $\mathbf{grad}U_t(\tilde{\theta}_t)$ around $\hat{\theta}_t$ yields

$$0 = -tc(\tilde{\theta}_t)\{\mathbf{grad}U_t(\hat{\theta}_t) + (\mathbf{H}U_t(\hat{\theta}_t))(\tilde{\theta}_t - \hat{\theta}_t)\}(1 + o_{\text{proba}}(1)) + \mathbf{grad}c(\tilde{\theta}_t).$$

In this equation, $\mathbf{grad}U_t(\hat{\theta}_t) = 0$ because $\hat{\theta}_t$ is the maximizer of $U_t(\cdot)$. Moreover, applying zero order Taylor's expansions for $c(\tilde{\theta}_t)$, $\mathbf{H}U_t(\hat{\theta}_t)$ and $\mathbf{grad}c(\tilde{\theta}_t)$ around $\theta$ yields

$$0 = -tc(\theta)(\mathbf{H}U_t(\theta))(\tilde{\theta}_t - \hat{\theta}_t)(1 + o_{\text{proba}}(1)) + \mathbf{grad}c(\theta)$$
$$= -tc(\theta)I_\theta(\tilde{\theta}_t - \hat{\theta}_t)(1 + o_{\text{proba}}(1)) + \mathbf{grad}c(\theta),$$

since $\lim_{t \to \infty} \mathbf{H}U_t(\theta) = I_\theta$ in probability. Then equation (4) follows.

## 3.4 Convergence of the CB–posterior distribution

The CB–posterior distribution $p_t(\cdot)$ is asymptotically equivalent to the density function of the Gaussian distribution $\mathcal{N}\left(\tilde{\theta}_t, (tI_\theta)^{-1}\right)$:

$$p_t(\alpha) \underset{t \to \infty}{\sim} \frac{1}{(2\pi)^{p/2}|(tI_\theta)^{-1}|^{1/2}} \exp\left(-\frac{1}{2}(\alpha - \tilde{\theta}_t)'(tI_\theta)(\alpha - \tilde{\theta}_t)\right). \qquad (5)$$

See the end of the section for the proof. This result allows us to figure out what is the CB–posterior distribution and how it can be used to make inference in the frequentist and Bayesian ways.

In the contrast theory, the distribution $\mathcal{N}\left(\tilde{\theta}_t, (tI)_\theta^{-1}\Gamma_\theta I_\theta^{-1}\right)$ is used to make frequentist inference about $\theta$: the point estimator is $\tilde{\theta}_t$, and confidence zones are provided based on the this normal distribution. Consequently, if the contrast is such that $I_\theta^{-1}\Gamma_\theta I_\theta^{-1} = I_\theta^{-1}$, then the CB–posterior distribution $p_t(\cdot)$ which approximates the density of $\mathcal{N}\left(\tilde{\theta}_t, (tI_\theta)^{-1}\right)$ can be directly used to make frequentist inference about $\theta$: the mode of $p_t(\cdot)$ is the point estimator, and confidence zones can be directly determined from $p_t(\cdot)$. This case is particularly interesting since the calculation of the limit matrices $I_\theta = \lim_{t\to\infty} \mathbf{H}U_t(\theta)$ and $\Gamma_\theta = \lim_{t\to\infty} V_\theta(\sqrt{t}\mathbf{grad}U_t(\theta))$ is not required.

Moreover, when the contrast which is considered satisfies $I_\theta^{-1}\Gamma_\theta I_\theta^{-1} = I_\theta^{-1}$, we propose to use the CB–posterior distribution $p_t(\cdot)$ to make inference in the Bayesian way, i.e. to use $p_t(\cdot)$ as a real posterior density. The motivation is based on the following analogy: when the contrast corresponding to the likelihood is employed (in this case, $I_\theta^{-1}\Gamma_\theta I_\theta^{-1} = I_\theta^{-1}$), then $p_t(\cdot)$ can be used (i) to make frequentist inference since it is an approximation of the limit distribution of the estimator (see above) and (ii) to make Bayesian inference since it is the classical posterior density. It has to be noted that, in the general case, the CB–posterior density $p_t(\cdot)$ does not coincide with the classical posterior density. It is a posterior density based on the information brought by the contrast under consideration.

If the contrast does not satisfy $I_\theta^{-1}\Gamma_\theta I_\theta^{-1} = I_\theta^{-1}$, then the CB–posterior distribution $p_t(\cdot)$ cannot be used to approximate the limit distribution of $\tilde{\theta}_t$ or to make Bayesian inference. However, $p_t(\cdot)$ can be used to estimate the matrix $I_\theta$, so avoiding the calculation of the second derivatives of the contrast. Indeed, one can see from (5) that an estimate of $I_\theta$ is the matrix $\Omega^{-1}/t$ where $\Omega$ is the variance matrix of the normal density function centered around $\tilde{\theta}_t$ and fitted to $p_t(\cdot)$. If $\theta$ is real, $I_\theta$ can be more simply estimated by $2\pi p_t(\tilde{\theta}_t)^2/t$ since equation (5) yields $p_t(\tilde{\theta}_t) \underset{t\to\infty}{\sim} (tI_\theta/2\pi)^{1/2}$. We have not found an equivalent way to easily estimate $\Gamma_\theta$ without analytical calculation of the second derivatives and without simulations.

Proof of (5). Let $\delta > 0$. For any $a$ such that $\sup_{1 \leq i \leq p} |a_i| < t^\delta$, a third order Taylor's expansion yields

$$\log p_t(\tilde{\theta}_t + a/\sqrt{t}) - \log p_t(\tilde{\theta}_t) = -\sqrt{t}a'\mathbf{grad}U_t(\tilde{\theta}_t) - \frac{1}{2}a'I_\theta a + o_{\text{proba}}(t^{2\delta} + t^{3\delta - 1/2}).$$

Given that $\mathbf{grad}U_t(\hat{\theta}_t) = 0$ (definition of the classical minimum contrast estimator $\hat{\theta}_t$) and that $\tilde{\theta}_t - \hat{\theta}_t = o_{\text{proba}}(t^{-1+\delta})\mathbf{1}_p$ (see eq. (4)), the previous equation becomes

$$\log p_t(\tilde{\theta}_t + a/\sqrt{t}) - \log p_t(\tilde{\theta}_t) = -\frac{1}{2}a'I_\theta a + o_{\text{proba}}(t^{2\delta} + t^{3\delta - 1/2}).$$

Ensuring that $\delta < 1/2$ (and not only $\delta > 0$), then

$$\log p_t(\tilde{\theta}_t + a/\sqrt{t}) - \log p_t(\tilde{\theta}_t) = -\frac{1}{2}a'I_\theta a + o_{\text{proba}}(t^{2\delta})$$
$$= -\frac{1}{2}a'I_\theta a \{1 + o_{\text{proba}}(1)\}.$$

Let us introduce $g_t : a \mapsto t^{-p/2} p_t(\tilde{\theta}_t + a/\sqrt{t})$ defined over $\mathbb{R}^p$. This density function satisfies, from the previous result,

$$g_t(a) \underset{t\to\infty}{\sim} t^{-p/2} p_t(\tilde{\theta}_t) \exp\left(-\frac{1}{2}a'I_\theta a\right).$$

Since $g_t(\cdot)$ is a density function and given the form of the right-hand-side term of this equation, $g_t(\cdot)$ is equivalent to the density function of the normal law with variance matrix $I_\theta^{-1}$. Equation (5) is then obtained with the change of variable $\alpha = \tilde{\theta}_t + a/\sqrt{t}$.

## 3.5 Summary: making inference with the CB–posterior distribution

For any contrast, a point estimator of $\theta$ is at the mode of the CB–posterior distribution $p_t(\cdot)$. Moreover, if $I_\theta^{-1}\Gamma_\theta I_\theta^{-1} = I_\theta^{-1}$, then $p_t(\cdot)$ can be used to make inference in the Bayesian way or in the frequentist way. Otherwise, $p_t(\cdot)$ can be used to estimate the limit matrix $I_\theta$.

It has to be noted that building a contrast such that $I_\theta^{-1}\Gamma_\theta I_\theta^{-1} = I_\theta^{-1}$ is particularly interesting since the calculation of $I_\theta$ and $\Gamma_\theta$ is avoided. However, we will see below that it is not always possible.

# 4 Applications in spatial statistics

## 4.1 Least-square estimation of a variogram range

This simulated case illustrates the application of the method for a real parameter. Here, the CB–posterior distribution cannot be directly used to make inference but can be used for estimating $I_\theta$.

We built a data set by simulating a centered Gaussian random field whose variogram is $\gamma_\theta(r) = 1 - \exp(-\theta r)$ with $\theta = 1$; $\theta$ is the inverse of the range parameter. The field was simulated over a $n \times n$ square grid ($n = 20$) with inter-node distance one. Figure 1 (left) shows the simulated random field. The sample variogram $\hat{\gamma}(h)$ was estimated for every possible inter-points distance $h$ less than the half diagonal of the grid; let $\mathcal{H}$ denote the set of these distances.

For the estimation of $\theta$, we chose a uniform prior density over $[0, 4]$ (horizontal dotted line in Fig. 1, right) and we used the least-square contrast introduced in section 2.2 (see eq. (1)). The CB–posterior density is shown in Figure 1 (right, dotted curve). The MAP estimate is $\tilde{\theta}_t = 1.34$ (vertical line).

Estimation uncertainty was assessed by estimating the limit variance of $\tilde{\theta}_t$ which is $\Gamma_\theta/(nI_\theta)^2$. The term $\Gamma_\theta = \lim_{t\to\infty} V_\theta(\sqrt{t}\mathbf{grad}U_t(\theta))$ ($t = n^2$ here) was estimated based on Monte-Carlo simulations: 1000 Gaussian random fields were simulated under $\tilde{\theta}_t$; for each simulation the sample variogram $\{\hat{\gamma}(h) : h \in \mathcal{H}\}$ was computed, and the first derivative of the contrast in $\tilde{\theta}_t$, i.e. $-\sum_{h\in\mathcal{H}} he^{-\tilde{\theta}_t h}\{\hat{\gamma}(h) - (1 - e^{-\tilde{\theta}_t h})\}$, was calculated; the sample variance of the derivatives multiplied by $n^2$ gave the estimate 1.97 for $\Gamma_\theta$.

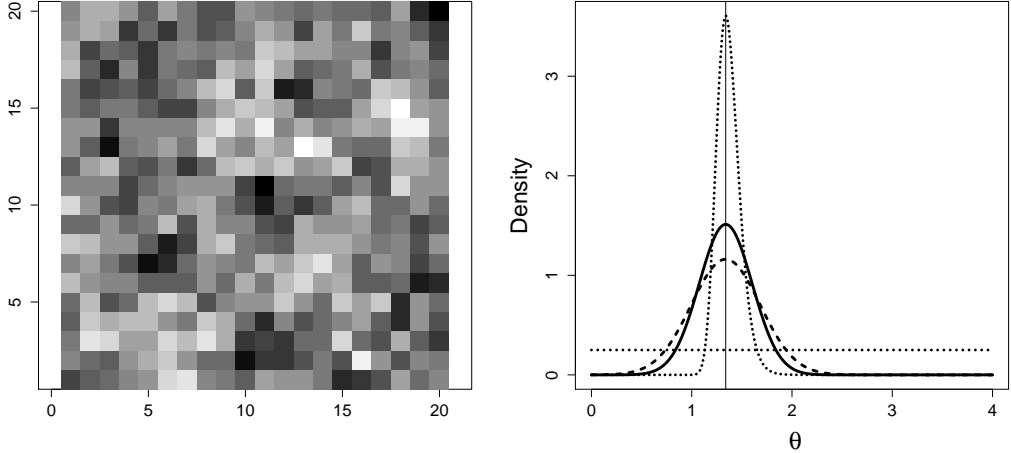

Figure 1: Left: realization of a centered Gaussian random field with exponential variogram parameterized by $\theta = 1$, over a $20 \times 20$ square-grid. Right: prior density (horizontal dotted line), contrast-based posterior density (dotted curve), density function of the limit distribution $\mathcal{N}(\tilde{\theta}_t, \Gamma_\theta/(nI_\theta)^2)$ (continuous and dashed lines when the estimate of the limit variance is based on simulations and when it is based on the posterior distribution), and MAP estimator (vertical line).

The term $I_\theta = \lim_{t \to \infty} \mathbf{H} U_t(\theta)$ was estimated in two ways: with the estimator $2\pi p_t(\tilde{\theta}_t)^2/t$ as suggested in section 3.4 and with Monte-Carlo simulations. In the former way, the estimate of $I_\theta$ is 0.20. The second way was carried out as follows: for each of the 1000 simulated Gaussian fields mentioned above, the second derivative of the contrast in $\tilde{\theta}_t$, i.e. $\sum_{h \in \mathcal{H}} h^2 e^{-\tilde{\theta}_t h}[e^{-\tilde{\theta}_t h} - \{\hat{\gamma}(h) - (1 - e^{-\tilde{\theta}_t h})\}]$, was computed; then, the sample mean of these derivatives gave the estimate 0.27 for $I_\theta$.

Thus, the estimate of the limit variance $\Gamma_\theta/(nI_\theta)^2$ of $\tilde{\theta}_t$ is 0.07 when $I_\theta$ is assessed by simulations and 0.12 when $I_\theta$ is computed from the CB–posterior distribution. The density function of the limit distribution $\mathcal{N}(\tilde{\theta}_t, \Gamma_\theta/(nI_\theta)^2)$ is drawn in Figure 1 (right). The continuous and dashed lines show this density when the estimate of the limit variance is 0.07 and 0.12, respectively. The true value $\theta = 1$ belongs to the 95%-confidence interval whatever the estimate of the limit variance is. We see how the two versions of the limit density are different from the CB–posterior density.

To assess the efficiency of the method, the coverage rate of the 95%-confidence interval was measured by applying the estimation procedure to 1000 simulated fields. The coverage rate is 94.6% when the estimate of $I_\theta$ is based on Monte-Carlo simulations and 94.7% when the estimate of $I_\theta$ comes from the contrast-based posterior density.

## 4.2 Pseudo-likelihood estimation of a Markovian spatial model

This simulated case illustrates the application of the method for a bivariate parameter. Here, the CB–posterior distribution is close from the limit distribution of the estimator. Here also, this posterior distribution cannot be directly used to make inference but can be used for estimating $I_\theta$.

We built a data set by simulating the spatial Markov field with two states, 0 and 1, specified in section 2.2. The field was simulated on a $n \times n$ square grid $\mathcal{I}$ ($n = 20$). Figure 2 (left) shows a simulation of this field for $\theta_1 = 0$ and $\theta_2 = 0.3$. To estimate $\theta_1$ and $\theta_2$, we applied the estimation method proposed in this article by using the pseudo-likelihood contrast introduced in section 2.2 (see eq. (2)) and a uniform prior density over $[-1.5, 1.5]^2$. The CB–posterior density is shown in Figure 2 (center). The MAP estimate is $\tilde{\theta}_t = (-0.21, 0.38)$.

For providing the limit distribution $\mathcal{N}(\tilde{\theta}_t, I_\theta^{-1} \Gamma_\theta I_\theta^{-1}/n^2)$ of the estimator, matrices $\Gamma_\theta$ and $I_\theta$ must be estimated. We computed the gradient and the Hessian of the contrast for $N = 1000$ Markov fields simulated under $\tilde{\theta}_t$, and we used the sample variance of the gradients for estimating $\Gamma_\theta$ and the sample mean of the Hessians for estimating $\Gamma_\theta$. The estimate of the limit variance matrix $I_\theta^{-1} \Gamma_\theta I_\theta^{-1}/n^2$

was finally
$$\begin{pmatrix} 0.14 & -0.055 \\ -0.055 & 0.022 \end{pmatrix}.$$

Almost the same matrix was obtained when we estimated $I_\theta$ by fitting a normal density to the CB–posterior density as suggested in section 3.4. Figure 2 (right) shows the limit density function of the estimator together with the 95%-confidence zone. We can see that the true parameter belongs to this zone. Moreover, Figure 2 shows the limit density is quite close from the posterior density. The pseudo-likelihood which accounts for short-distance interactions certainly brings almost the same information than the likelihood brings. It has however to be noted that this would not be the case if long-distance interactions had been introduced in the spatial Markov model.

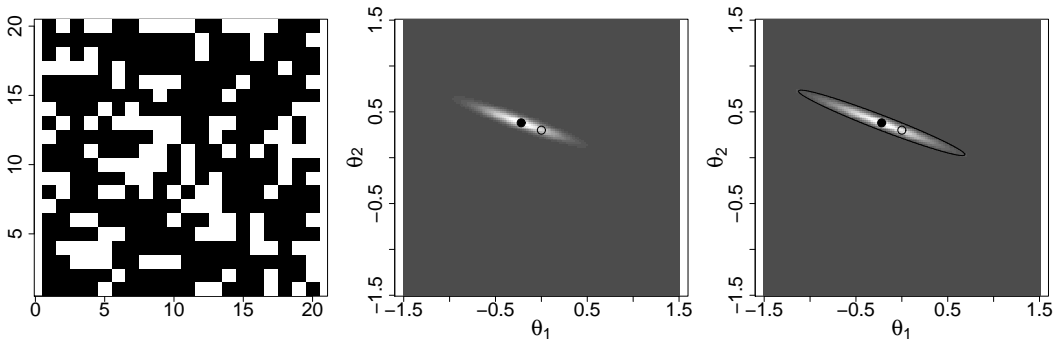

Figure 2: Left: realization of a Markovian spatial process with two states over a 20×20 grid. Center: contrast-based posterior density. Right: limit density $\mathcal{N}(\tilde{\theta}_t, I_\theta^{-1}\Gamma_\theta I_\theta^{-1}/n^2)$. On the center and right panels, the MAP estimate and the true parameter are drawn with a black dot and a circle, respectively. On the right panel, the continuous line circumscribes the 95%-confidence zone.

### 4.3 Estimation of an autosimilar model using moments

This real study-case illustrates the application of the method for a bivariate parameter. Here, the CB–posterior distribution can be directly use to make inference.

In this section we aims to build and fit a model for soil roughness. Soil roughness plays an important role in the distribution of rain water into infiltration, pond and streaming. It also modifies reflectance properties of soils used to estimate soil moisture with remote detection for example. An experiment was carried out to measure soil roughness at a small scale. Soil heights were measured every 2mm along 1.18m-transects in a cultivated field (Bertuzzi et al., 1995). Figure 3 (top) shows the distributions of heights for two among twelve sampled transects. These distributions were obtained after subtraction of the trend estimated with a kernel smoothing. The mean height computed from the 12 transects is 7.6mm, the maximum is 22.9mm. Several models have been proposed to describe soil surface. For instance, in Boolean models and autosimilar models (Bertuzzi et al., 1995; Goulard and Chadœuf, 1994; Lantuéjoul, 2002, chap. 14), basic random elements (e.g. cylinder) are drawn from a given law and the soil surface is the maximum height in the former model and the summed height in the latter model.

Here, we aim to estimate the parameters of an autosimilar model based on random cylinders, each cylinder having same height and radius. For any $x \in \mathbb{R}^2$ and $r > 0$, let $f(x, r) = r 1_{\{||x|| < r\}}$ be the function describing the cylinder which is centered in $x$ and whose radius and height are equal to $r$. In addition, let $(X, R)$ be a marked Poisson point process defined over $\mathbb{R}^2 \times \mathbb{R}^{+*}$ with intensity function $\mu(x, r) = \alpha \exp\{-\beta r\}$. The random surface $Y$ representing the soil surface is defined by

$$Y_M = \sum_{(x,r) \in (X,R)} f(x - M, r).$$

For such a process, it is difficult to calculate the joint distribution of the heights whereas the moments can easily be written. The parameter vector $\theta = (\alpha, \beta)$ has two components and we propose to estimate it using the first two moments: $\hat{\mu}_A = (\frac{1}{\nu(A)} \int_A Y_M dM, \frac{1}{\nu(A)} \int_A Y_M^2 dM)$, where $A$ is the set of the sampled transects and $\nu(A)$ is its measure.

If border effects are neglected, the expected value of $\hat{\mu}_A$ is

$$E(\hat{\mu}_A) = \left(6\pi\frac{\alpha}{\beta^4}, 36\pi^2\frac{\alpha^2}{\beta^8} + 24\pi\frac{\alpha}{\beta^5}\right).$$

Moreover, the variance matrix of $\hat{\mu}_A$ satisfies

$$\nu(A)\text{var}(\hat{\mu}_A) \to V,$$

where the components of $V$ are

$$V_{11} = 5!\frac{16}{3}\frac{\alpha}{\beta^6}$$
$$V_{12} = 6!\frac{16}{3}\frac{\alpha}{\beta^7} + (5!)64\pi\frac{\alpha^2}{\beta^{10}}$$
$$V_{22} = 7!\frac{16}{3}\frac{\alpha}{\beta^8} + \{(6!)128\pi + (10!)32\kappa\}\frac{\alpha^2}{\beta^{11}} + (3!)(5!)128\pi^2\frac{\alpha^3}{\beta^{14}},$$

with $\kappa = \int_0^1 \int_0^1 (\arccos(u) - u\sqrt{1-u^2})(\arccos(v) - v\sqrt{1-v^2})\frac{(uv)^5}{(u+v)^{11}}dudv$.

The estimation method is applied by using a uniform prior over $[1,100]\times[1,5]$ and a contrast based on the weighted least squares of the first two moments:

$$U_A(\theta) = (\hat{\mu}_A - E(\hat{\mu}_A))'V^{-1}(\hat{\mu}_A - E(\hat{\mu}_A))/2.$$

For this contrast, the matrices $I_\theta$ and $\Gamma_\theta$ are equal and their component $(i,j)$ is

$$\frac{\partial E(\hat{\mu}_A)'}{\partial \theta_i}V^{-1}\frac{\partial E(\hat{\mu}_A)}{\partial \theta_j}.$$

Consequently, $I_\theta^{-1}\Gamma_\theta I_\theta^{-1} = I_\theta^{-1}$) and the CB–posterior density can be used as an approximation of the limit density of the MAP estimator $\tilde{\theta}_A$ or as a posterior distribution of the parameter $\theta$ (see section 3.4). Figure 3 (bottom) shows the joint CB–posterior distribution and the marginals. The MAP estimate of $\theta$ is $\tilde{\theta}_A = (46.6, 3.28)$. Marginal 95%-confidence intervals of $\alpha$ and $\beta$ are [36.1,58.5] and [3.07,3.48], respectively.

## 5 Discussion

We have studied a method of estimation exploiting a contrast-based posterior distribution (CBPD). This method includes the classical likelihood-based proce-dures (MLE and Bayesian estimation), but has been mainly developed to circumvent difficulties encountered with the likelihood by generalizing the Bayesian formula of the posterior distribution, so extending the proposal of Lin (2006). The CBPD can be used to make frequentist inference and, in specific situations, Bayesian inference. In case of frequentist inference, the use of the CBPD allows the reduction of analytical calculations usually required to compute the limit variance matrix of the estimator. In this article, the method has been applied to spatial data sets, but can be applied to other cases where likelihood-based procedures are not appropriate.

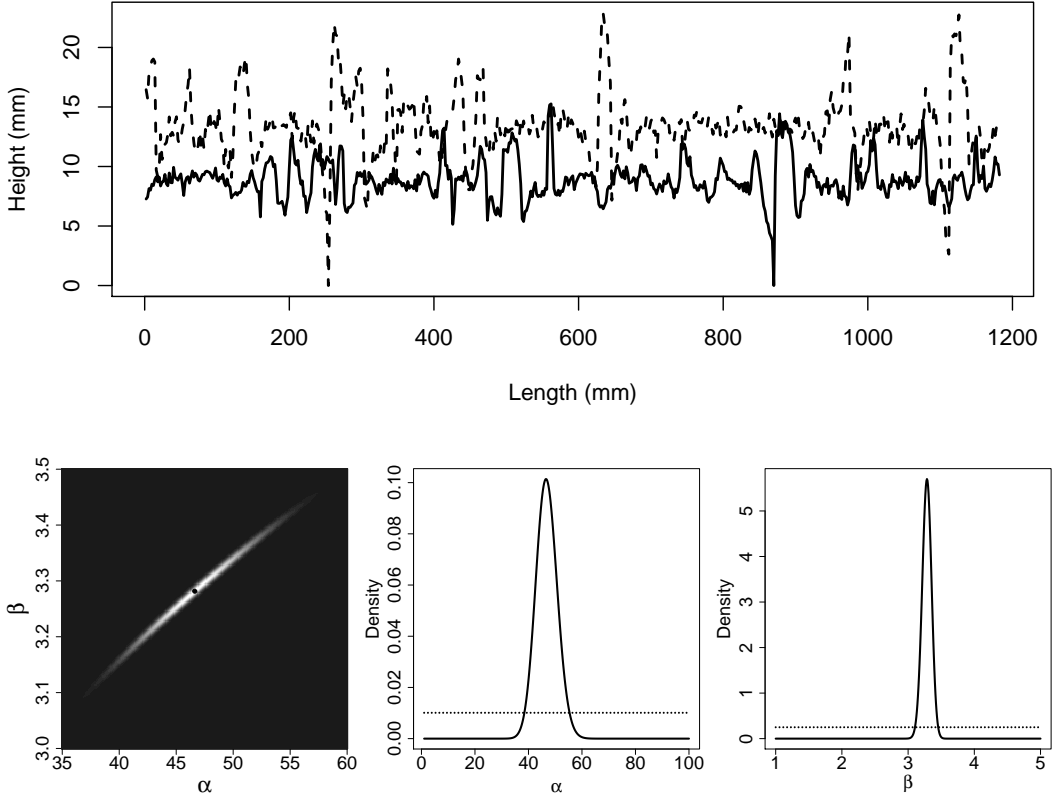

Figure 3: Top: distribution of heights for two transects (heights were corrected by kernel smoothing for subtracting the trend). Bottom left: contrast-based posterior density for $(\alpha, \beta)$; the MAP estimate is at the black dot. Bottom center and right: contrast-based posterior marginal densities for $\alpha$ et $\beta$ (continuous lines) and prior marginal densities (dashed lines).

In the frequentist viewpoint, the CBPD can be used to provide a point estimator (the posterior mode) and the limit distribution of this estimator. The limit distribution is directly approximated by the CBPD if the variance of the gradient vector of the contrast is equal to the inverse of the limit Hessian matrix of the contrast (i.e. $\Gamma_\theta = I_\theta^{-1}$; see the third application). In this case, it is not required to calculate and estimate the variance matrix of the estimator. In other cases, the limit distribution is not directly available, but the Hessian matrix of the contrast can be easily estimated from the CBPD and, consequently, the calculation of the second derivatives of the contrast is avoided (see the first two applications). It has to be noted that using Bayesian calculation to make frequentist estimation has been proposed in the literature (Robert and Hwang, 1996; Robert and Titterington, 1998; Jacquier et al., 2007), but the proposals were restricted to maximum likelihood estimation.

In the Bayesian viewpoint, the CBPD can be used as a classical posterior distribution when $\Gamma_\theta = I_\theta^{-1}$, as in the third application. It has however to be noted that the CBPD does not always coincide with the classical posterior distribution. The CBPD has to be viewed as a posterior distribution based on the information brought by the contrast which is used.

Even if the proposed procedure has advantages, it also faces two classical limits: the choice of the prior distribution (or the penalization function in the frequentist viewpoint) which can influence the posterior inference, and the choice of the contrast. Regarding the former limit, we refer to Clarke and Gustafson (1998) and Rootzén and Olsson (2006) for example. Regarding the choice of the contrast, we have two comments. The first comment concerns the possibility to build a contrast such that $\Gamma_\theta = I_\theta^{-1}$ (case where our method is the most interesting). It was possible in the real case-study because we could provide the analytical form for the variance matrix of the sample moments. However, it was not possible in the two simulated case-studies. Indeed, for the estimation of the range parameter, we should have modeled the variance of the variogram. However, such a practice is not common in geostatistics when the field is not assumed to be Gaussian. For the estimation of the spatial Markov model, the spatial dependences make impossible to get a transformed pseudo-likelihood such that $\Gamma_\theta = I_\theta^{-1}$; it has to be noted that the problem of dependence can be circumvented with coding techniques (Besag, 1975) but, with such techniques, a part of the information is lost. This leads us to our second comment about the information brought by contrasts. We see that in the real case-study the two estimators are strongly correlated. We could have tried to use another contrast to avoid correlation. For example, together with the sample mean, we could have used the covariance at a given distance instead of the variance to get two moments which are less correlated. However, the calculation of the expected value and the variance-covariance of these moments is much more tricky. Thus, to be able to derive analytical expressions and apply the method as it is presented, the choice of the contrast is limited. Nevertheless, simulations could be used to circumvent this difficulty as in approximate Bayesian computation (Beaumont et al., 2002). This could be an interesting extension of the estimation method proposed in this paper.